\documentclass[11pt]{amsart}
\usepackage[]{amsmath,amssymb,fullpage}

\begin{document}
\title {Generalized Bunce--Deddens algebras}
\author [S. Orfanos]  {Stefanos Orfanos}
\address  {Department of Mathematical Sciences, University of Cincinnati, Cincinnati, OH, 45221}
\email {stefanos.orfanos@uc.edu}

%Abstract
\begin{abstract}
We define a broad class of crossed product C*-algebras of the form $C(\tilde{G})\rtimes G$, where $G$ is a discrete countable amenable residually finite group, and $\tilde{G}$ is a profinite completion of $G$. We show that they are unital separable simple nuclear quasidiagonal C*-algebras, of real rank zero, stable rank one, with comparability of projections and with a unique trace.
\end{abstract}

\maketitle

\theoremstyle{plain}
\newtheorem {thm} {Theorem} 
\newtheorem {lem} [thm]{Lemma} 
\newtheorem {cor} [thm]{Corollary}
\newtheorem {pro} [thm]{Proposition}

\theoremstyle{definition}
\newtheorem {dfn} {Definition}
\newtheorem {que} {Question}
\newtheorem* {ex} {Example}

\theoremstyle{remark}
\newtheorem* {rem} {Remark}
\newtheorem* {akn} {Aknowledgments}
\newtheorem* {pr} {Proof}

% Introduction
\section*{Introduction}
Our object of study will be a family of crossed products, which we call the \emph{generalized Bunce--Deddens algebras}, because their construction generalizes a well-known construction of the classical Bunce--Deddens algebras, with $\mathbb{Z}$ being replaced by any discrete countable amenable residually finite group $G$, and the odometer action appropriately generalized. Such actions have been considered in the literature, in topological dynamics and quite recently in von Neumann algebra theory, but the corresponding crossed product C*-algebras have not been explicitly described. We should mention here that the generalized Bunce--Deddens algebras appearing in this paper are different from the ones introduced by D. Kribs in \cite{krib1}. 

The generalized Bunce--Deddens algebras turn out to have many desirable properties, being simple nuclear quasidiagonal, and having real rank zero, stable rank one, a unique trace and comparability of projections. The ultimate goal of classifying such algebras is not realized in this paper; however, it is conjectured that they will have tracial rank zero or finite decomposition rank. Another open problem is to compute their ordered K--theory (possibly, their Elliott invariant), which may be within reach for certain groups.

% First section
\section{Amenable residually finite groups and profinite completions}

We start with the definition of amenability. Throughout this paper, $G$ will be assumed to be discrete and countable.
\begin{dfn} 
A group $G$ is \emph{amenable} if there is a sequence $e\in F_1\subseteq F_2\subseteq\dotsc\subseteq F_n\dotsc$ of finite sets of $G$ such that:
\[\bigcup_{n\ge 1}F_n = G \mbox{ and }\lim_{n\to\infty}\frac{|F_n\triangle F_ns|}{|F_n|}=0\mbox{ for all }s\in G\]
\end{dfn}

Every amenable group has also left F\o lner sequences, i.e., sequences $(F_n)_n$ that also exhaust the group $G$ and satisfy \[\lim_{n\to\infty}\frac{|F_n\triangle sF_n|}{|F_n|}=0\mbox{ for all }s\in G\]
as well as sequences that are both left and right F\o lner. Unless noted otherwise, a \emph{F\o lner sequence} will be a right F\o lner sequence from now on.
\begin{dfn}
The group $G$ is \emph{residually finite} if it has a separating family of finite index normal subgroups.
\end{dfn}

In other words, for every finite set $F$ in $G$, there is a normal subgroup $L$ of finite index in $G$, such that the quotient map $G\to G/L$ is injective, when restricted to $F$.

A \emph{tiling} of $G$ is a decomposition $G=KL$, where $K$ is called the \emph{tile} and $L$ is the set of \emph{tiling centers}, so that every $x\in G$ is uniquely written as a product of an element in $K$ and an element in $L$. Tilings are an important feature of amenable residually finite groups, as we will see below.

\begin{lem}
Suppose $G$ is an amenable and residually finite group, and assume a separating nested sequence of finite index subgroups $L_n$ of $G$ is given. Then there exists a F\o lner sequence $(F_n)_n$ and a sequence of finite subsets $K_n\supset F_n$ such that $G$ has a tiling of the form $G=K_nL_n$ for all $n\ge 1$. 
\end{lem}

The proof is based on the previous definitions. We will need a stronger result concerning the asymptotic behavior of tiles for $G$. The following theorem, essentially due to B. Weiss from \cite{weis1}, constructs \emph{F\o lner tiles} for every separating nested sequence of finite index normal subgroups of $G$. The statement below comes from \cite{desc1}.

\begin{thm}[Weiss]
For every discrete countable amenable residually finite group $G$ and every separating nested sequence of finite index normal subgroups $(L_n)_n$, there is a sequence of sets $(K_n)_n$ that is left and right F\o lner, such that $G=K_nL_n$ is a tiling for all $n\ge 1$. 
\end{thm}
 
Proofs of this result can be found in \cite{desc1} or \cite{weis1}.

In what follows, we refine the above result to get left/right F\o lner sets that tile, not only the group $G$, but also subsequent tiles. First, we recall an equivalent formulation of F\o lner's characterization of amenability.

\begin{dfn}
The \emph{$S$-boundary} of a set $K$, denoted $\partial_SK$, is the set $SK\cap SK^c = \{ x\in G : S^{-1}x\cap K \ne \emptyset \mbox{ and } S^{-1}x\cap K^c \ne \emptyset \}$. 
\end{dfn}

\begin{lem}The left F\o lner condition for $K_n$ can be written: \[\frac{|\partial_SK_n|}{|K_n|}\to 0 \mbox { as } n\to\infty \mbox{ for any fixed finite }S\subset G.\]
\end{lem}

Indeed, $SK_n = (SK_n \cap SK_n^c) \cup (SK_n \setminus SK_n^c) = (SK_n \cap SK_n^c) \cup (\displaystyle \bigcap _{s\in S}sK_n)$, and both $\displaystyle \frac{|SK_n|}{|K_n|}$, $\displaystyle \frac{|\displaystyle \cap _{s\in S}sK_n|}{|K_n|}$ converge to 1 by Theorem 16.16 of \cite{pier1}.

The improved tiling result is the following:

\begin{thm}
If $G$ is a discrete countable amenable residually finite group $G$ and $(L_n)_n$ is a separating nested sequence of finite index normal subgroups of $G$, then there exist finite subsets $K_n$ of $G$ and integers $n_1 <n_2 < \dotsm$ such that:
\begin{enumerate}
\item $(K_n)_n$ is a left and right F\o lner sequence, and $G=K_nL_n$ is a tiling for all $n\ge 1$
\item $K_{n_{k+1}}$ is the disjoint union of right translates of $K_{n_k}$, for all $k=1, 2, \dotsc$.
\end{enumerate} 
\end{thm} 

\begin{pr}
The first property comes from Weiss's theorem. To prove the second property, let $\epsilon>0$, a finite set $S\in G$, and an integer $n_k>0$ be given. We claim that there is $n_{k+1}>n_k$ such that the F\o lner tile $K_{n_{k+1}}$ and the set $K_{n_{k+1}}'= K_{n_k}(K_{n_{k+1}}\cap L_{n_k})$ satisfy the following condition:
\[\frac{|K_{n_{k+1}} \triangle K_{n_{k+1}}'|}{|K_{n_{k+1}}'|} < \frac{\epsilon}{3}\]
To prove the claim, note that $|K_{n_{k+1}}'|=|K_{n_k}||K_{n_{k+1}}\cap L_{n_k}| = |K_{n_k}|[L_{n_k}:L_{n_{k+1}}] = |K_{n_{k+1}}|$. We also observe that $G=K_{n_k}L_{n_k}$ is partitioned into three sets: 
\begin{itemize}
\item $G_1$ = union of right translates of $K_{n_k}$ that are subsets of $K_{n_{k+1}}$ 
\item $G_2$ = union of right translates of $K_{n_k}$ that are subsets of $K_{n_{k+1}}^c$ 
\item $G_3$ = union of right translates of $K_{n_k}$ that intersect both $K_{n_{k+1}}$ and $K_{n_{k+1}}^c$ 
\end{itemize}
Then, $G_1\subset K_{n_{k+1}}\cap K_{n_{k+1}}' \subset K_{n_{k+1}}\cup K_{n_{k+1}}'\subset G_2^c$ and $G_2^c \setminus G_1 = G_3$. Consequently,
\[K_{n_{k+1}}\triangle K_{n_{k+1}}' \subset G_3 =\displaystyle \bigcup\{K_{n_k}l :l\in L_{n_k}\mbox{ with } K_{n_k}l\cap K_{n_{k+1}} \ne \emptyset \mbox{ and                                            } K_{n_k}l\cap K_{n_{k+1}}^c \ne \emptyset \}\]
\[ \subset K_{n_k} \partial_{K_{n_k}^{-1}}K_{n_{k+1}} = K_{n_k} (K_{n_k}^{-1}K_{n_{k+1}} \cap K_{n_k}^{-1}K_{n_{k+1}}^c) \subset \partial_{K_{n_k}K_{n_k}^{-1}}K_{n_{k+1}}\] 
Hence, by the previous lemma, there exists $n_{k+1}>n_k$ such that the claim is true. Enlarge $n_{k+1}$ if necessary, so that 
\[ |K_{n_{k+1}}\triangle K_{n_{k+1}}s|<\frac{\epsilon}{3}|K_{n_{k+1}}| \mbox{ and } |K_{n_{k+1}}\triangle sK_{n_{k+1}}|<\frac{\epsilon}{3}|K_{n_{k+1}}| \mbox{, for all }s\in S\]
Then we claim that $K'_{n_{k+1}}$ satisfies:
\[ |K'_{n_{k+1}}\triangle K'_{n_{k+1}}s|<\epsilon |K'_{n_{k+1}}| \mbox{ and } |K'_{n_{k+1}}\triangle sK'_{n_{k+1}}|<\epsilon |K'_{n_{k+1}}| \mbox{, for all }s\in S\]
We show the left F\o lner condition below; the same argument works for the right F\o lner condition too. First observe that
\[K_{n_{k+1}}'\triangle sK_{n_{k+1}}' \subset (K_{n_{k+1}}'\triangle K_{n_{k+1}})\cup (K_{n_{k+1}}\triangle sK_{n_{k+1}})\cup (sK_{n_{k+1}}\triangle sK_{n_{k+1}}')\] to get \[|K_{n_{k+1}}'\triangle sK_{n_{k+1}}'| \le |K_{n_{k+1}}'\triangle K_{n_{k+1}}|+ |K_{n_{k+1}}\triangle sK_{n_{k+1}}|+ |sK_{n_{k+1}}\triangle sK_{n_{k+1}}'|\]\[ <\frac{\epsilon}{3} |K_{n_{k+1}}'|+\frac{\epsilon}{3} |K_{n_{k+1}}'|+\frac{\epsilon}{3} |K_{n_{k+1}}'| = \epsilon |K_{n_{k+1}}'| \mbox{ for all } s\in S\]
Therefore, $K_{n_{k+1}}'$ is a left/right F\o lner tile which is tiled by $K_{n_k}$, and thus we can replace $K_{n_{k+1}}$ by it.\qed
\end{pr}

We now discuss profinite completions of groups. Let $I$ be a directed set. We define an \emph{inverse system} $(G_i,\phi_{ij})$ as a net of topological groups $G_i$, $i\in I$ and continuous homomorphisms $\phi_{ij}:G_j\to G_i$ such that $\phi_{ii}$ is the identity morphism, for all $i\in I$ and $\phi_{ij}\circ \phi_{jk} = \phi_{ik}$ for all $i\le j\le k$. In the category of topological groups, inverse systems have always a unique limit. The \emph{inverse limit} of the inverse system $(G_i, \phi_{ij})$, denoted by $\displaystyle \lim_{\gets}G_i$, is the subgroup of the product $\displaystyle \prod_{i\in I}G_i$ consisting of sequences $(x_i)_i$ that satisfy $\phi_{ij}(x_j)=x_i$ for all $i\le j$. 

In particular, for $G$ a residually finite group, we fix a decreasing sequence of finite index normal subgroups $L_n$ that separates points of $G$. The groups $G_n$ in the definition of the inverse system are the finite quotients $G/L_n$, and the homomorphisms $\phi_{nm}:G/L_m\to G/L_n$ are given by $\phi_{nm}(xL_m)=xL_n$ for $n\le m$. Then the \emph{profinite completion} of $G$ with respect to these subgroups, denoted by $\tilde{G}$, is the inverse limit of the finite quotients $G/L_n$, that is, the subgroup of $\displaystyle \prod_{n\ge 1}G/L_n$ consisting of sequences $(x_nL_n)_n$ such that $x_mL_n= \phi(x_mL_m)=x_nL_n$ whenever $n\le m$. 

We denote by $\pi_n$ the canonical projections onto $G/L_n$, $n\ge 1$ and we formulate below some well-known properties of the profinite completion of $G$ (cf. \cite{wils1}). 

\begin{pro}The profinite completion has the following properties:
\begin{enumerate}
\item $\tilde{G}$ is a non-empty totally disconnected compact Hausdorff group.
\item The sets $\pi_n^{-1}(\{xL_n\})$, $xL_n\in G/L_n$ and $n\ge 1$, form a base of compact and open sets for the topology on $\tilde{G}$.
\end{enumerate}
\end{pro}

\begin{pr}
The first assertion follows from standard general topology arguments. To prove the second assertion, note that the sets $\pi_n^{-1}(\{xL_n\})$ are compact and open for all $xL_n\in G/L_n$ and $n\ge 1$. Let $U$ be open in $\tilde{G}$ and $(x_nL_n)_n\in U$. Then, there are integers $n_1<\dotsm <n_m$ and open sets $U_{n_1}, \dotsc,U_{n_m}$ in the respective quotients, such that $\pi_{n_1}^{-1}(U_{n_1})\cap\dotsm\cap\pi_{n_m}^{-1}(U_{n_m})\subset U$. In particular, $\pi_{n_1}^{-1}(\{x_{n_1}L_{n_1}\})\cap\dotsm\cap\pi_{n_m}^{-1}(\{x_{n_m}L_{n_m}\})\subset U$. But $\pi_{n_1}^{-1}(\{x_{n_1}L_{n_1}\})\supset\dotsm\supset\pi_{n_m}^{-1}(\{x_{n_m}L_{n_m}\})$, hence the conclusion. \qed
\end{pr}

A useful corollary is that every compact set in $\tilde{G}$ is a finite union of the above-mentioned base sets. Also, $G$ embeds as a dense subgroup into $\tilde{G}$, and similarly, $L_n$ is dense in $\ker \pi _n$ for all $n\ge 1$. For ease of notation, we will denote $\ker \pi_n$ by $\tilde{L_n}$ from now on. Finally, we observe that $\pi^{-1}_n(\{ xL_n\} ) = x\tilde{L_n}$ for every $x\in G$ (after the identification of $G$ as a subgroup of $\tilde{G}$).  

% Second section
\section{Construction of the Generalized Bunce--Deddens algebras}

We assume that the reader is familiar with the construction of crossed products. A good reference is \cite{will1}. In our case, the reduced crossed product $A\rtimes _{\alpha ,r}G$ coincides with the full crossed product $A\rtimes _{\alpha}G$, since $G$ is an amenable group.

Let $q=(q_n)_n$ be a sequence of positive integers such that $q_{n+1}$ is divisible by $q_n$ for all $n\ge 1$. The usual way to define the \emph{Bunce--Deddens algebra of type $q$} is to consider the sequence of finite groups
\[\mathbb{Z}/q_1\mathbb{Z}\to \mathbb{Z}/q_2\mathbb{Z}\to\dotsm\to\mathbb{Z}/q_n\mathbb{Z}\to\dotsm\]
where each term acts on $C(\mathbb{T})$ by rotations, and take the inductive limit of the corresponding crossed products. Since $C(\mathbb{T})\rtimes\mathbb{Z}/q_n\mathbb{Z}$ is isomorphic to the algebra of $q_n\times q_n$ matrices with entries in $C(\mathbb{T})$, we can also consider the resulting algebra as the inductive limit of such matrix algebras. 

A second approach is to start with the limit of the inverse sequence

\[\mathbb{Z}/q_1\mathbb{Z}\gets \mathbb{Z}/q_2\mathbb{Z}\gets \dotsm\gets\mathbb{Z}/q_n\mathbb{Z}\gets\dotsm\]
which is just the profinite completion $\tilde{\mathbb{Z}}$ of $\mathbb{Z}$ with respect to these subgroups. Then, consider the integers as a subgroup of $\tilde{\mathbb{Z}}$ and define the action by addition. The crossed product $C(\tilde{\mathbb{Z}})\rtimes\mathbb{Z}$ is again the Bunce--Deddens algebra of type $q$, as seen by the fact that the action by rotations of the inductive limit of $\mathbb{Z}/q_n\mathbb{Z}$'s on $\mathbb{T}$ induces the action by addition of the group $\tilde{\mathbb{T}}$ ($\cong\mathbb{Z}$) of characters of $\mathbb{T}$ on the inverse limit of the $\widetilde{\mathbb{Z}/q_n\mathbb{Z}}$ ($\cong\mathbb{Z}/q_n\mathbb{Z}$).

We now generalize this second construction in the following way: Let $G$ be an amenable residually finite group with a sequence of nested finite index normal subgroups $L_n$ that separates points, and consider its action $\alpha$ by left multiplication on its profinite completion $\tilde{G}$ with respect to these subgroups. The resulting crossed products $C(\tilde{G})\rtimes _{\alpha} G$ are the \emph{generalized Bunce--Deddens algebras}. A few of their properties are straightforward to check. Since $\tilde{G}$ is compact, $C(\tilde{G})$ is unital, and together with $G$ being discrete, makes $C(\tilde{G})\rtimes G$ unital as well. Separability is also evident: $\tilde{G}$ is metrizable, therefore $C(\tilde{G})$ is separable, and since $G$ is countable, $C(\tilde{G})\rtimes G$ is clearly separable too.

Recall that an action of $G$ on a locally compact Hausdorff space $X$ is \emph{free} if every stabilizer $\{g\in G:gx=x\}$ is trivial, and \emph{minimal} if every orbit $\{gx : g\in G\}$ is dense in $X$. Since $G$ is amenable and its action on $\tilde{G}$ is free and minimal, the crossed product $C(\tilde{G})\rtimes G$ is a simple C*-algebra, and nuclearity is a consequence of the theorem below (refer to \cite{broz1} for proofs). 

\begin{thm}
A C*-algebra $A$ is nuclear, if and only if there exist contractive completely positive maps $\phi_n: A\to M_{k_n}(\mathbb{C})$ and $\psi_n:M_{k_n}(\mathbb{C})\to A$ such that their composition approximates the identity map in the point-norm topology. In particular, if $G$ is a discrete amenable group acting on a compact Hausdorff space $X$, then the crossed product $C(X)\rtimes G$ is nuclear.
\end{thm}

Therefore, we have concluded that:

\begin{cor}
The generalized Bunce--Deddens algebras are unital simple separable and nuclear.
\end{cor} 

Finally, we claim that the generalized Bunce--Deddens algebras are quasidiagonal. We start with the definition of quasidiagonality.
\begin{dfn}
A linear operator $T$ on a separable Hilbert space $\mathcal{H}$ is \emph{quasidiagonal} if there exists a sequence of finite rank self-adjoint orthogonal projections $P_n$ in $\mathcal{B(H)}$ satisfying:
\begin{enumerate}
\item $ P_n\to I_{\mathcal{H}}$ as $n\to\infty$, and
\item $\|[T,P_n]\| \to 0$ as $n\to\infty$ 
\end{enumerate}
A separable set of operators $\mathcal{A}$ is \emph{quasidiagonal} if every operator $T$ in a set of dense linear span in $\mathcal{A}$ is quasidiagonal with respect to the same sequence $(P_n)_n$. An abstract C*-algebra $A$ is \emph{quasidiagonal} if it has a faithful representation to a quasidiagonal set of operators.
\end{dfn}
We use the following theorem from \cite{orfa1}.
\begin{thm}
Let $G$ be a discrete countable amenable and residually finite group with a sequence of F\o lner sets $F_n$ and tilings of the form $G=K_nL_n$ with $F_n\subset K_n$ for all $n\ge 1$. Let $A$ be a unital separable C*-algebra and let $\alpha :G \to Aut A$ be a homomorphism such that \[\lim _{n\to\infty} \left[\max_{l\in L_n \cap K_nK_n^{-1}F_n}\|\alpha (l)a - a\|\right] = 0\]for all $a\in A$. Assume, moreover, that $A$ is quasidiagonal. Then $A\rtimes _\alpha G$ is also quasidiagonal.
\end{thm}

Observe that due to Proposition~5 and the remarks after it, every function $f$ in $C(\tilde{G})$ can be approximated within $\epsilon$ by the sum of constant functions supported on the compact open sets $x\tilde{L_n}$, where $x\in K_n$ and with $n$ depending on $\epsilon$ and the function $f$. Moreover, recall that the action of any element $l\in L_n$ on $x\tilde{L_n}=\tilde{L_n}x$ leaves it invariant, since $L_n$ is embedded in $\tilde{L_n}$ for all $n\ge 1$. Hence, 
\[\alpha (l)\chi _{x\tilde{L_n}} = \chi _{x\tilde{L_n}} \mbox{ for any } l\in L_n \]and by the triangle inequality, 
\[\max_{l\in L_n}\|\alpha (l)f - f\|<2\epsilon \mbox{ for $n$ large enough}\]
It follows that the action is almost periodic, as defined in the statement of Theorem~8, and thus we have the following:

\begin{thm}
Every generalized Bunce--Deddens algebra is quasidiagonal.
\end{thm}

% Third section
\section{Almost AF groupoids and further properties}

In this section we introduce groupoids and focus especially on transformation groups as such. The connection with crossed products is the following: The groupoid algebra of a transformation group $X\rtimes G$ is canonically isomorphic to the crossed product $C(X)\rtimes G$. 
\begin{dfn}
A \emph{groupoid} is a set $\mathcal{G}$ with an associative product defined on the subset of $\mathcal{G}\times \mathcal{G}$ consisting of \emph{composable pairs}, and an inverse defined everywhere. The inverse satisfies $(x^{-1})^{-1}=x$ and every element makes a composable pair with its inverse (in either order), but also $x^{-1}x$ need not equal $xx^{-1}$. The pair $(x,y)$ is composable if and only if $x^{-1}x = yy^{-1}$, in which case both $x^{-1}xy=y$ and $xyy^{-1}=x$ are true. The set of elements of the form $x^{-1}x$ for $x\in \mathcal{G}$ is the \emph{unit space} of $\mathcal{G}$, denoted by $\mathcal{G}^0$. The element $s(x)=xx^{-1}$ is called the \emph{source} of $x\in \mathcal{G}$, and $r(x)=x^{-1}x$ is its \emph{range}.
\end{dfn}

We focus on $\mathcal{G}=\tilde{G}\rtimes G$ from now on. The action is by left multiplication, and the elements of $\tilde{G}\rtimes G$ are of the form $(x,g)$, with $x\in \tilde{G}$ and $g\in G$. The product is defined on pairs $((x,g),(y,h))$ such that $y=gx$, by the formula $(x,g)(y,h)=(x,hg)$, and the inverse $(x,g)^{-1} = (gx, g^{-1})$. $\tilde {G}\rtimes G$ is thus a groupoid, with unit space isomorphic to $\tilde{G}$. Think of $(x,g)$ as an arrow from $x\in \tilde{G}$ to $gx\in \tilde{G}$. The fact that all arrows are defined uniquely by their endpoints is a consequence of the action being free, and makes $\tilde{G}\rtimes G$ a \emph{principal} groupoid. Also, it is a \emph{Cantor} groupoid, which is defined as a second countable locally compact Hausdorff etale groupoid, whose unit space is the Cantor set, equipped with a Haar system of counting measures. Indeed, for transformation groups $X\rtimes G$ it is sufficient that $X$ is homeomorphic to the Cantor set and that $G$ is discrete and countable (cf. \cite{phil1}). Moreover, an open subgroupoid of a Cantor groupoid with the same unit space is itself a Cantor groupoid. 

We continue with the definion of an AF groupoid from \cite{gips1}

\begin{dfn}[Giordano--Putnam--Skau]
A Cantor groupoid $\mathcal{G}$ is called \emph{approximately finite} (AF) if it is the increasing union of a sequence of compact open principal Cantor subgroupoids, with each of them containing $\mathcal{G}^0$.
\end{dfn}

We show that $\mathcal{G}=\tilde{G}\rtimes G$ contains an AF groupoid, by constructing a nested sequence of compact open subgroupoids. Based on the tiling $G=K_nL_n$ of the group $G$, and the fact that $G$ acts freely on $\tilde{G}$, we get a tiling
\[\tilde {G} = \displaystyle \bigcup_{xL_n\in G/L_n} \pi^{-1}_n(\{xL_n\}) =  K_n\tilde{L_n}\] 
of its profinite completion. We then define the following subsets of $\mathcal{G}$:
\[\mathcal{G}_n = \{(gx, g_1g^{-1}): x\in \tilde{L_n}, \, g,g_1\in K_n\}\]
It is easy to check that $\mathcal{G}_n$ are subgroupoids of $\mathcal{G}$ with the same unit space. They are compact and open in $\mathcal{G}$ as a consequence of $\tilde{L_n}$ being such. In addition, Theorem~4 allows us to obtain a subsequence of these subgroupoids which is nested, i.e., to find integers $n_1 < n_2 < \dotsm$ such that $\mathcal{G}_{n_k}\subset \mathcal{G}_{n_{k+1}}$ for all $k\ge 1$. Indeed, if $(gx, g_1g^{-1})\in \mathcal{G}_{n_k}$, then write $gx=hy\in K_{n_{k+1}}\tilde{L_{n_{k+1}}}$. Since $K_{n_{k+1}}$ has a partition of the form $\bigcup_{l\in L} K_{n_k}l$ for a suitable finite set $L\subset L_{n_k}$, we have $K_{n_{k+1}}y = \bigcup_{l\in L} K_{n_k}ly$, hence $gx \in K_{n_k}l_1y$ for some $l_1 \in L$, which gives $x=l_1y$ (since the action is free). Then, $h_1=g_1l_1\in K_{n_{k+1}}$ and thus $(gx, g_1g^{-1}) = (hy,h_1h^{-1}) \in \mathcal{G}_{n_{k+1}}$.

The groupoid algebra of an AF groupoid is an AF algebra, and conversely, to any AF algebra we can associate a unique AF groupoid with groupoid algebra the AF algebra we started with (cf. \cite{rena1}). Observe however that for elements of infinite order in a group $G$, the corresponding elements in $C^*(G)$ cannot have finite spectrum. Therefore, $\mathcal{G} = X\rtimes G$ is not an AF groupoid, as long as $G$ is not a locally finite group. Yet, it may happen that $\mathcal{G}$ is not `much bigger' than an open AF subgroupoid, in a sense that N. C. Phillips made precise in \cite{phil1}. The following definition applies when the groupoid algebra of $\mathcal{G}$ is simple. A \emph{graph} is a subset of $\mathcal{G}$ for which the restrictions of the source and range maps are injective.

\begin{dfn}[Phillips]
Assume $\mathcal{G}$ is a Cantor groupoid such that $C_r^*(\mathcal{G})$ is simple. Then $\mathcal{G}$ is \emph{almost AF} if it contains an open subgroupoid $\mathcal{G}_{AF}$ with the same unit space and the following property: for every compact subset $C$ of $\mathcal{G}\setminus \mathcal{G}_{AF}$ and every $m\ge 1$, there exist compact graphs $C_1, \dotsc ,C_m$ in $\mathcal{G}_{AF}$ with source $s(C_i)=s(C)$ for $i=1,\dotsc ,m$ and disjoint ranges.
\end{dfn}  

The crucial step is to show that $\mathcal{G}=\tilde{G}\rtimes G$ is an almost AF groupoid. A more general result of this sort appears in \cite{phil1}; however Phillips has to assume finite generation (which we do not need).

\begin{thm}
The groupoid $\mathcal{G}=\tilde{G}\rtimes G$ is almost AF for every discrete countable amenable residually finite group $G$ and every profinite completion $\tilde{G}$ associated to a separating nested sequence of finite index normal subgroups $L_n$ of $G$.
\end{thm}

\begin{pr}
We are inspired from the proof of Theorem~6.9 in \cite{phil1}. We will use the already established notation. We will also identify the unit space of $\mathcal{G}$ with $\tilde{G}$ to simplify notation. Set $\mathcal{G}_{AF} = \bigcup_{k\ge 1} \mathcal{G}_{n_k}$ which is open in $\mathcal{G}$. The goal is to verify the last condition of Definition~7. To that end, consider a compact set $C\subset \mathcal{G} \setminus \mathcal{G}_{AF}$ and an integer $m\ge 1$. Define
\[ S =\{g\in G:(\tilde{G}\times \{g\})\cap C\ne \emptyset\}\]
which is a finite set because $C$ is compact. 

By Theorem~4, there exists $n\in \{n_1, n_2, \dotsc\}$ so that $|K_n\triangle sK_n|< \displaystyle \frac{1}{|S|m}|K_n|$ for all $s\in S$. For $(y,s)\in C$, $y\in \tilde{G}=K_n\tilde{L_n}$, so $y=gx$ for some $g\in K_n$ and $x\in \tilde{L_n}$. However, $(y,s)\notin \mathcal{G}_n$, hence $(y,s)\ne (gx, g_1g^{-1})$ for any $g_1\in K_n$. It follows that $g\in K_n\setminus s^{-1}K_n$. If $K=\displaystyle \bigcup_{s\in S} (K_n\setminus s^{-1}K_n)$, then we have $s(C)\subset \displaystyle \bigcup_{g\in K} g\tilde{L_n}$ with
\[|K|\le \sum_{s\in S} |K_n\setminus s^{-1}K_n|\le |S||K_n\triangle sK_n|<\frac{1}{m}|K_n|\]
Therefore, there exist $m$ injective functions $\sigma_1, \dotsc , \sigma_m:K\to K_n$ with disjoint ranges. As a result, the compact sets $C_i= \displaystyle \bigcup_{g\in K} [(s(C)\cap g\tilde{L_n})\times \{\sigma_i(g)g^{-1}\}]$, $i=1, \dotsc , m$, satisfy 
\[s(C_i)=\displaystyle \bigcup_{g\in K} (s(C)\cap g\tilde{L_n})=s(C)\] 
and their ranges 
\[r(C_i) = \bigcup_{g\in K} [\sigma_i(g)g^{-1}s(C)\cap \sigma_i(g)\tilde{L_n}] \subset \bigcup_{g\in K} \sigma_i(g)\tilde{L_n}\]
are disjoint by definition of the functions $\sigma_i$. 

Moreover, if $(gx,\sigma_i(g)g^{-1})$ and $(hy,\sigma_i(h)h^{-1})$ are two elements of $C_i$, then \[s(gx,\sigma_i(g)g^{-1}) = gx = hy = s(hy,\sigma_i(h)h^{-1}) \mbox{ forces } g=h \mbox{ and } x=y\]and the same is true for the ranges. Hence $C_i$ is a graph for all $i=1,\dotsc , m$. Finally, $C_i\subset \mathcal{G}_n$ because the functions $\sigma _i$ map into $K_n$ for every $i=1,\dotsc , m$. \qed
\end{pr} 

The significance of almost AF groupoids comes from Phillips's theorem below. But first, let us recall the following notions:

\begin{dfn}
A (unital) C*-algebra has \emph{real rank zero} if the self-adjoint invertible elements are dense in the set of all self-adjoints (Brown--Pedersen, \cite{brpe1}). A C*-algebra has \emph{stable rank one} if the invertible elements are dense in the algebra (Rieffel, \cite{rief1}). 
\end{dfn}

\begin{thm}[Phillips]
Let $\mathcal{G}$ be an almost AF groupoid, and assume $C_r^*(\mathcal{G})$ is simple. Then $C_r^*(\mathcal{G})$ has real rank zero and stable rank one. Moreover, $C_r^*(\mathcal{G})$ has comparability of projections: for $p,q$ projections in $\mathcal{M}_\infty(C_r^*(\mathcal{G}))$ with $\tau (p) <\tau (q)$ for all normalized traces $\tau$ on $C_r^*(\mathcal{G})$, $p$ is Murray--von Neumann equivalent to a subprojection of $q$. Finally, there is a bijection between normalized traces on $C_r^*(\mathcal{G})$ and invariant Borel probability measures on $\mathcal{G}^0$.
\end{thm}

A Borel measure $\mu$ on $\mathcal{G}^0$ is \emph{invariant} if for every $f\in C_c(\mathcal{G})$, the following is true:
\[
\int_{ \mathcal{G}^0} \left(\sum_{g\in\mathcal{G}:s(g)=x} f(g)\right) d\mu (x) = \int_{ \mathcal{G}^0} \left(\sum_{g\in\mathcal{G}:r(g)=x} f(g)\right) d\mu (x)\]

The correspondence between traces and invariant measures on the unit space is obtained for more general groupoids and is more explicit in \cite{phil1}.

We are now able to prove:

\begin{thm}
Every generalized Bunce--Deddens algebra has real rank zero, stable rank one, comparability of projections, and a unique trace.
\end{thm}

\begin{pr}
We combine Theorems~10 and~11 to get real rank zero, stable rank one and comparability of projections. The unique trace is obtained as follows: Observe that
\[ S =\{g\in G:(\tilde{G}\times \{g\})\cap \mbox{supp}f\ne \emptyset\}\]
is finite, and hence, that the integrals in the equation that gives the invariance of a Borel measure on $\mathcal{G}^0$ are:
\[\sum_{s\in S}\int_{\tilde{G}} f(x,s) d\mu (x) \mbox {, and } \sum_{s\in S}\int_{\tilde{G}} f(s^{-1}x,s) d\mu (x)\]
We now see that the two integrals are equal if and only if the measure $\mu$ is $G$-invariant. However, we can show that the normalized Haar measure is the only $G$-invariant probability measure on $\tilde{G}$. Indeed, if $\mu$ is any $G$-invariant measure, $f\in C_c(\tilde{G})$, and $(y_n)_n$ is a sequence of elements of $G$ converging to $y\in \tilde{G}$, then 
\[ \int_{\tilde{G}} f(x) d\mu (x) = \int_{\tilde{G}} f(y_n^{-1}x) d\mu (x) \longrightarrow \int_{\tilde{G}} f(y^{-1}x) d\mu (x) \mbox{, as }n\to \infty\]Therefore $\mu$ is also $\tilde{G}$-invariant, and we are done.
\qed
\end{pr}  

% Aknowledgements
\begin{akn}
The results presented here are part of the author's doctoral thesis at Purdue University. I~am indebted to Marius Dadarlat for tons of advice and support and I~would like to thank Chris Phillips for many helpful discussions. Last but not least, many thanks to Larry Brown for comments on an earlier draft.
\end{akn}

% References

\end{document}